\documentclass[a4paper,12pt]{amsart}
\usepackage{amsfonts}
\usepackage{amssymb}
\usepackage{ifthen}
\usepackage{graphicx}

\newtheorem{thm}{Theorem}
\newtheorem{cor}{Corollary}
\newtheorem{lem}{Lemma}
\newtheorem{prop}[equation]{Proposition}
\newtheorem{examp}{Example}

\newtheorem{conj}[equation]{Conjecture}
\newtheorem{rem}{Remark}
\theoremstyle{definition}
\newtheorem{defn}{Definition}

\newtheorem{prob}[equation]{Problem}
\newtheorem{ques}[equation]{Question}

\newcounter {own}
\def\theown {\thesection       .\arabic{own}}

\newenvironment{pf}[1][]{%
 \vskip 3mm
 \noindent
 \ifthenelse{\equal{#1}{}}%
  {{\slshape Proof. }}%
  {{\slshape #1.} }%
 }%
{\qed\bigskip}

\newcounter{alphabet}
\newcounter{tmp}
\newenvironment{Thm}[1][]{\refstepcounter{alphabet}%
\bigskip%
\noindent%
{\bf Theorem \Alph{alphabet}}%
\ifthenelse{\equal{#1}{}}{}{ (#1)}%
{\bf .} \itshape}{\vskip 8pt}

\newenvironment{Lem}[1][]{\refstepcounter{alphabet}%
\bigskip%
\noindent%
{\bf Lemma \Alph{alphabet}}%
{\bf .} \itshape}{\vskip 8pt}

\def\be{\begin{equation}}
\def\ee{\end{equation}}

\newcommand{\bee}{\begin{enumerate}}
\newcommand{\eee}{\end{enumerate}}

\newcommand{\blem}{\begin{lem}}
\newcommand{\elem}{\end{lem}}
\newcommand{\bthm}{\begin{thm}}
\newcommand{\ethm}{\end{thm}}
\newcommand{\bcor}{\begin{cor}}
\newcommand{\ecor}{\end{cor}}
\newcommand{\beg}{\begin{examp}}
\newcommand{\eeg}{\end{examp}}
\newcommand{\begs}{\begin{examples}}
\newcommand{\eegs}{\end{examples}}
\newcommand{\bdefe}{\begin{defn}}
\newcommand{\edefe}{\end{defn}}
\newcommand{\bprob}{\begin{prob}}
\newcommand{\eprob}{\end{prob}}
\newcommand{\bques}{\begin{ques}}
\newcommand{\eques}{\end{ques}}
\newcommand{\bei}{\begin{itemize}}
\newcommand{\eei}{\end{itemize}}

\newcommand{\bcon}{\begin{conj}}
\newcommand{\econ}{\end{conj}}
\newcommand{\bcons}{\begin{conjs}}
\newcommand{\econs}{\end{conjs}}
\newcommand{\bprop}{\begin{prop}}
\newcommand{\eprop}{\end{prop}}
\newcommand{\br}{\begin{rem}}
\newcommand{\er}{\end{rem}}
\newcommand{\brs}{\begin{rems}}
\newcommand{\ers}{\end{rems}}
\newcommand{\bo}{\begin{obser}}
\newcommand{\eo}{\end{obser}}
\newcommand{\bos}{\begin{obsers}}
\newcommand{\eos}{\end{obsers}}
\newcommand{\bpf}{\begin{pf}}
\newcommand{\epf}{\end{pf}}
\newcommand{\ba}{\begin{array}}
\newcommand{\ea}{\end{array}}
\newcommand{\beq}{\begin{eqnarray}}
\newcommand{\beqq}{\begin{eqnarray*}}
\newcommand{\eeq}{\end{eqnarray}}
\newcommand{\eeqq}{\end{eqnarray*}}

\begin{document}
\bibliographystyle{amsplain}
\title [] {Landau's theorem for polyharmonic mappings}

\author{J. Chen}
\address{J. Chen, Department of Mathematics,
Hunan Normal University, Changsha, Hunan 410081, People's Republic
of China.} \email{jiaolongchen@sina.com.cn}

\author{ A. Rasila }
\address{A. Rasila, Department of Mathematics,
Hunan Normal University, Changsha, Hunan 410081, People's Republic
of China, and
Department of Mathematics and Systems Analysis, Aalto University, P. O. Box 11100, FI-00076 Aalto,
 Finland.} \email{antti.rasila@iki.fi}

\author{X. Wang${}^{~\mathbf{*}}$}
\address{X. Wang, Department of Mathematics,
Hunan Normal University, Changsha, Hunan 410081, People's Republic
of China.} \email{xtwang@hunnu.edu.cn}

\keywords{polyharmonic  mapping, Landau's theorem.\\
${}^{\mathbf{*}}$ Corresponding author}

\begin{abstract}
In this paper, we first investigate coefficient estimates for  bounded polyharmonic mappings in the unit disk $\mathbb{D}$. Then, we obtain two versions of Landau's theorem for polyharmonics mapping $F$, and for the mappings of the type $L(F)$, where $L$ is the differential operator of Abdulhadi, Abu Muhanna and Khuri. Examples and numerical estimates are given.
\end{abstract}

\thanks{The research was partly supported by
NSF of China (No. 11071063) and  Hunan Provincial Innovation
Foundation for Postgraduate (No. 125000-4242). }

\maketitle \pagestyle{myheadings} \markboth{J. Chen,  A. Rasila and X. Wang}
{Landau's theorem for planar polyharmonic mappings
}

\section{Introduction and preliminaries}\label{csw-sec1}

A $2p$ $(p\geq1)$ times continuously differentiable
complex-valued mapping $F=u+iv$ in a domain $D \subset \mathbb{C}$ is said to be {\it
polyharmonic} (or {\it $p$-harmonic}) if $F$ satisfies the
polyharmonic equation $\Delta^{p}F =\Delta(\Delta^{p-1}F)= 0$,
where $\Delta:=\Delta^{1}$ is the usual complex Laplacian operator
$$\Delta=4\frac{\partial^{2}}{\partial z\partial \overline{z}}:=
\frac{\partial^{2}}{\partial x^{2}}+\frac{\partial^{2}}{\partial y^{2}}. $$
Obviously, for $p=1$ (resp. $p=2$), we obtain the usual class of harmonic (resp. biharmonic) mappings.
The reader is referred to \cite{sh2011, CRW, qi} for discussion on biharmonic and polyharmonic mappings, and \cite{cl, du} for the general theory of harmonic mappings.

The biharmonic equation arises from certain problems in physics, in particular,
fluid dynamics and elasticity theory. This equation also has important applications in engineering and
biology (cf. \cite{ha, kh, la}). The classes of biharmonic and polyharmonic mappings can be understood as  natural generalizations of the class of harmonic mappings. However, the investigation of these mappings in the geometric function theory has  been started only recently
(cf. \cite{ab2008, ab2005, ab2006, sh2009, sh2010, CRW,CRW2,CW}).

Let $\mathbb{D}_{r} =\{z:|z|<r\}$ $(r >0$), and denote by $\mathbb{D}$ the unit disk $\mathbb{D}_{1}$. The classical theorem of Landau states that if $f$ is analytic in $\mathbb{D}$, with $f(0)=f'(0)-1=0$
and $|f(z)|<M$ for some $M\geq1$, then $f$ is univalent in the disk $\mathbb{D}_{\rho}$ with
$$\rho=\frac{1}{M+\sqrt{M^{2}-1}}>0.$$ In addition, the range $f(\mathbb{D}_{\rho})$ contains a
disk of radius $M\rho^{2}$ (cf. \cite{ch}).  Recently, there has
 been a number of articles dealing with Landau's theorem for planar harmonic mappings, see, for
example, \cite{ch, sh2010, sh, do, gr}, and for biharmonic mappings,
see \cite{ab2008, sh2009, sh2010, li}.

In this paper, we study coefficient estimates for certain classes of polyharmonic mappings. Our goal is to prove Landau type results for these mappings. The main results are Theorems \ref{thm1.1}, \ref{thm2.1} and \ref{thm3.1}, which are presented in Section \ref{csw-sec2}. Theorem \ref{thm1.1} is a generalization of \cite[Lemma $2.3$]{li2012}, Theorem \ref{thm2.1} is a generalization of \cite[Theorem 2]{sh2011} and Theorem \ref{thm3.1} is a
generalization of \cite[Theorem 1.1]{sh2009}.

\section{Preliminaries}

If $D$ is a simply connected domain, then it is easy to see that a mapping
$F$ is polyharmonic if and only if $F$  has the following representation:
$$F(z)=\sum_{k=1}^{p}|z|^{2(k-1)}G_{k}(z),$$ where each $G_{k}$ is a
complex-valued harmonic mapping in $D$ for $k\in \{1,\cdots,p\}$
(cf. \cite{sh2011, qi}). It is known (cf. \cite{cl, du}) that the
mappings $G_{k}$ can be expressed as the form
\begin{center}
$G_{k} = h_{k} + \overline{g_{k}}$
\end{center}
for $k\in \{1,\cdots,p\}$, where all $h_{k}$ and  $g_{k}$ are
analytic in $D$. In this paper, we
consider the polyharmonic mappings in
$\mathbb{D}$.

For a polyharmonic mapping $F$ in $\mathbb{D}$, we use the following standard notations:
$$\lambda_{F}(z)=\min_{0\leq\theta\leq2\pi}|F_{z}(z)
+e^{-2i\theta}F_{\overline{z}}(z)|=\big||F_{z}(z)|-|F_{\overline{z}}(z)| \big|,$$
 $$\Lambda_{F}(z)=\max_{0\leq\theta\leq2\pi}|F_{z}(z)
 +e^{-2i\theta}F_{\overline{z}}(z)|= |F_{z}(z)|+|F_{\overline{z}}(z)| .$$
Thus, for a sense-preserving polyharmonic mapping $F$, one has $J_{F}(z)=\Lambda_{F}(z)\lambda_{F}(z).$

We now recall some lemmas, which gives us some coefficient estimates.

\begin{Lem}{\rm \bf \cite[Lemma $2.1$]{li2009}}\label{lemA}
Suppose that $f=h+\overline{g}$
is a harmonic mapping of $\mathbb{D}$ with $h(z)=\sum_{n=1}^{\infty}a_{n}z^{n}$ and
$g(z)=\sum_{n=1}^{\infty}b_{n}z^{n}$ for $z\in \mathbb{D}$.
If $J_{F}(0)=1$ and $|f(z)|< M$, then
$$|a_{n}|,\;|b_{n}|\leq \sqrt{M^{2}-1},\;\;n=2,3,\cdots,$$
$$|a_{n}|+|b_{n}|\leq \sqrt{2M^{2}-2},\;\;n=2,3,\cdots,$$
and
$$\lambda_{F}(0)\geq\lambda(M)=\begin{cases}
\displaystyle\frac{\sqrt{2}}{\sqrt{M^{2}-1}+\sqrt{M^{2}+1}}\;&\;\text{if}\;
1\leq M\leq\frac{\pi}{2\sqrt[4]{2\pi^{2}-16}},\\
\displaystyle\frac{\pi}{4M} &\;\text{if}\;M>\frac{\pi}{2\sqrt[4]{2\pi^{2}-16}}.
\end{cases}
$$
\end{Lem}

\begin{Lem}{\rm \bf \cite[Lemma $1$]{sh2011}}\label{lemB}
Let $f=h+\overline{g}$
be a harmonic mapping of $\mathbb{D}$ such that $|f(z)|< M$ with $h(z)=\sum_{n=0}^{\infty}a_{n}z^{n}$ and
$g(z)=\sum_{n=1}^{\infty}b_{n}z^{n}$. Then $|a_{0}|\leq M$ and for any $n\geq1$,
\be\label{eq000000.2}|a_{n}|+|b_{n}|\leq \frac{4M}{\pi}.\ee
The estimate \eqref{eq000000.2} is sharp. The extremal functions are $f(z)\equiv M$ or
$$f_{n}(z)=\frac{2M\alpha}{\pi}\arg \left(\frac{1+\beta z^{n}}{1-\beta z^{n}}\right),$$
where $|\alpha|=|\beta|=1$.
\end{Lem}

\begin{Lem}{\rm \bf \cite[Lemma $2.3$]{li2012}}\label{lemC}
Suppose that $f=h+\overline{g}$
is a harmonic mapping of the unit disk $\mathbb{D}$ with $h(z)=\sum_{n=0}^{\infty}a_{n}z^{n}$ and
$g(z)=\sum_{n=1}^{\infty}b_{n}z^{n}$.

$(I)$ If $f$ satisfies $|f(z)|\leq M$ for all $z\in \mathbb{D}$, then
$$\sum_{n=2}^{\infty}(|a_{n}|^{2}+|b_{n}|^{2})\leq M^{2}-(|a_{1}|^{2}+|b_{1}|^{2}).$$

$(II)$ If $f$ satisfies $|f(z)|\leq M$ for all $z\in \mathbb{D}$ and $|J_{F}(0)|=1$, then
$$\left(\sum_{n=2}^{\infty}(|a_{n}|^{2}+|b_{n}|^{2})\right)^{\frac{1}{2}}
\leq \sqrt{M^{4}-1}\cdot\lambda_{F}(0).$$

$(III)$ If $f$ satisfies $|f(z)|\leq M$ for all $z\in \mathbb{D}$ and $\lambda_{F}(0)=1$, then
$$\left(\sum_{n=2}^{\infty}(|a_{n}|^{2}+|b_{n}|^{2})\right)^{\frac{1}{2}}
\leq \sqrt{2(M^{2}-1)}.$$
\end{Lem}

In \cite{ab2006}, the authors considered the following differential operator
$L$ defined on the class of complex-valued $C^{1}$ functions:
$$L=z\frac{\partial}{\partial z}-\overline{z}\frac{\partial}{\partial \overline{z}}.$$
Clearly, $L$ is a complex linear operator, and $L$ satisfies the usual product rule:
$$L(af+bg)=aL(f)+bL(g)\;\mbox{and}\;L(fg)=fL(g)+gL(f),$$
where $a$, $b$ are complex constants, and $f$, $g$ are $C^{1}$
functions. In addition, the operator $L$ has a number
of interesting properties. For instance, it is easy to see
that the operator $L$ preserves both harmonicity and biharmonicity.
Other basic properties of this operator are presented in \cite{ab2006}.

\begin{Thm}{\rm \bf \cite[Corollary $1(3)$]{ab2006}}\label{thmB} Let
$F$ be a univalent biharmonic function
in $\mathbb{D}$. If F is convex and $L(F)$ is univalent, then $L(F)$ is
starlike.
\end{Thm}

It is also important to recall that (see \cite[Corollary $1(3)$]{ab2006})
the operator $L(F)$ for polyharmonic mappings behaves
much like $zf'(z)$ for analytic functions. From Theorem D,
we see that it is meaningful to consider Landau's theorem for
polyharmonic mappings of the form $L(F)$, where $F$ belongs to the class of polyharmonic mappings.

\section{Main results}\label{csw-sec2}
We begin this section with some coefficient estimates for polyharmonic mappings.

\begin{thm}\label{thm1.1}
Suppose that $F$ is a polyharmonic mapping of the form
\be\label{eq1.1} F(z)=a_{0}+\sum_{k=1}^{p}|z|^{2(k-1)}
(h_{k}(z)+\overline{g_{k}(z)})=a_{0}+\sum_{k=1}^{p}|z|^{2(k-1)}
\sum_{n=1}^{\infty}(a_{n,k}z^{n}+\overline{b_{n,k}}\overline{z^{n}}),\ee
and all its non-zero coefficients $a_{n,k_{1}}$, $a_{n,k_{2}}$ and
$b_{n,k_{3}}$, $b_{n,k_{4}}$ satisfy the condition:
\be\label{eq1.2} \left|\arg \frac{a_{n,k_{1}}}{a_{n,k_{2}}}\right|\leq\frac{\pi}{2}\;
\mbox{and}\;\left |\arg \frac{b_{n,k_{3}}}{b_{n,k_{4}}}\right|\leq\frac{\pi}{2}.\ee

$(I)$ If $F$ satisfies $|F(z)|\leq M$ for all $z\in \mathbb{D}$, then \be\label{eq1.3}|a_{0}|^{2}+\sum_{k=1}^{p}\sum_{n=1}^{\infty}(|a_{n,k}|^{2}+|b_{n.k}|^{2})\leq M^{2},\ee
and $|a_{n,k}|+|b_{n,k}|\leq\sqrt{2}M$ for all coefficients.

$(II)$ If $F$ satisfies $|F(z)|\leq M$ for all $z\in \mathbb{D}$ and $F(0)=|J_{F}(0)|-1=0$, then \be\label{eq1.4}\sqrt{\sum_{k=1}^{p}\sum_{n=1}^{\infty}(|a_{n,k}|+|b_{n.k}|)^{2}
-(|a_{1,1}|+|b_{1,1}|)^{2}}\leq \sqrt{M^{4}-1}\cdot \lambda_{F}(0),\ee
and $|a_{n,k}|+|b_{n,k}|\leq T_{1}(M):=\min \left\{\sqrt{2M^{2}-2},
\;\sqrt{M^{4}-1}\cdot \lambda_{F}(0)\right\}$
for all $(n,k)\not=(1,1)$, where
\be\label{eq1.7}\lambda_{F}(0)\geq\lambda_{0}(M):=\frac{\sqrt{2}}{\sqrt{M^{2}-1}+\sqrt{M^{2}+1}}.\ee

$(III)$ If $F$ satisfies $|F(z)|\leq M$ for all $z\in \mathbb{D}$ and $F(0)=\lambda_{F}(0)-1=0$, then \be\label{eq1.6}\sqrt{\sum_{k=1}^{p}\sum_{n=1}^{\infty}(|a_{n,k}|+|b_{n.k}|)^{2}
-(|a_{1,1}|+|b_{1,1}|)^{2}}\leq \sqrt{2M^{2}-2},\ee
and $|a_{n,k}|+|b_{n,k}|\leq \sqrt{2M^{2}-2}$
for all $(n,k)\not=(1,1)$.
\end{thm}

\bpf We first prove the inequality \eqref{eq1.3}. A standard argument using Parseval's identity and
the hypothesis that $|F(re^{i \theta})|\leq M$, gives
\begin{eqnarray*}
&&\frac{1}{2\pi}\int^{2\pi}_{0}|F(re^{i \theta})|^{2}d \theta\\
&=&|a_{0}|^{2}+\sum_{k=1}^{p}r^{4(k-1)}\sum_{n=1}^{\infty}(|a_{n,k}|^{2}+|b_{n,k}|^{2})r^{2n}\\
&&+2\mbox{Re}\sum_{1\leq k_{1}<k_{2}\leq p}r^{2(k_{1}+k_{2}-2)}\sum_{n=1}^{\infty} (a_{n,k_{1}}\overline{a_{n,k_{2}}}+b_{n,k_{1}}\overline{b_{n,k_{2}}})r^{2n}\\
&\leq&M^{2}.
\end{eqnarray*}
Since $F $ satisfies the condition \eqref{eq1.2}, it follows that $\mbox{Re}(a_{n,k_{1}}\overline{a_{n,k_{2}}}+b_{n,k_{1}}\overline{b_{n,k_{2}}})\geq0$
for any $k_{1}$, $k_{2}\in \{1,\cdots,p\}$.
Let $r\longrightarrow 1^{-}$. We obtain that
$$|a_{0}|^{2}+\sum_{k=1}^{p}\sum_{n=1}^{\infty}(|a_{n,k}|^{2}+|b_{n.k}|^{2})\leq M^{2},$$
which proves that the inequality \eqref{eq1.3} holds.
By Cauchy's inequality, for any coefficient $a_{n,k}$, $b_{n,k}$,
we get $$(|a_{n,k}|+|b_{n,k}|)^{2}\leq2(|a_{n,k}|^{2}+|b_{n,k}|^{2})\leq 2M^{2}.$$
That is $|a_{n,k}|+|b_{n,k}|\leq\sqrt{2}M$.

Then, we prove the inequality \eqref{eq1.4}. We have from \eqref{eq1.3} that
\begin{align}\begin{split}\label{eq1.5}
&\left(\sum_{k=1}^{p}\sum_{n=2}^{\infty}(|a_{n,k}|+|b_{n.k}|)^{2}
+\sum_{k=2}^{p}(|a_{1,k}|+|b_{1.k}|)^{2}\right)^{\frac{1}{2}}\\
\leq&\left(\sum_{k=1}^{p}\sum_{n=2}^{\infty}2(|a_{n,k}|^{2}+|b_{n.k}|^{2})
+\sum_{k=2}^{p}2(|a_{1,k}|^{2}+|b_{1.k}|^{2})\right)^{\frac{1}{2}}\\
\leq&\big(2M^{2}-2(|a_{1,1}|^{2}+|b_{1,1}|^{2})\big)^{\frac{1}{2}}.\\
\end{split}\end{align}

It follows from \eqref{eq1.5} that
$|a_{n,k}|+|b_{n.k}|\leq \sqrt{2M^{2}-2}$ for $(n,k)\not=(1,1)$.
Then if $J_{F}(0)=1$, we get $|a_{1,1}|=\sqrt{1+|b_{1,1}|^{2}}\geq1$.
 Thus we deduce
from \eqref{eq1.3} that $|b_{1,1}|\leq\sqrt{\frac{M^{2}-1}{2}}$, $M\geq1$.
Hence $$\sqrt{2M^{2}-2(|a_{1,1}|^{2}+|b_{1,1}|^{2})}=
\sqrt{2M^{2}-2(1+2|b_{1,1}|^{2})}\Big(|b_{1,1}|+\sqrt{1+|b_{1,1}|^{2}}\Big)\cdot \lambda_{F}(0).$$
As in \cite{li2012}, by a straightforward calculation, we get
$$\sqrt{2M^{2}-2(1+2|b_{1,1}|^{2})}\Big(|b_{1,1}|+\sqrt{1+|b_{1,1}|^{2}}\Big)\leq\sqrt{M^{4}-1},$$
which proves that the inequality \eqref{eq1.4} holds,
$|a_{n,k}|+|b_{n.k}|\leq \sqrt{M^{4}-1}\cdot \lambda_{F}(0)$ for $(n,k)\not=(1,1)$ and $$\lambda_{F}(0)=\big||a_{1,1}|-|b_{1,1}|\big|=\frac{1}{\sqrt{1+|b_{1,1}|^{2}}+|b_{1,1}|}
\geq\frac{\sqrt{2}}{\sqrt{M^{2}-1}+\sqrt{M^{2}+1}}.$$

If $J_{F}(0)=-1$, we get $|b_{1,1}|=\sqrt{1+|a_{1,1}|^{2}}\geq1$.
 Thus we deduce from \eqref{eq1.3} that $|a_{1,1}|\leq\sqrt{\frac{M^{2}-1}{2}}$, $M\geq1$.
Hence $$\sqrt{2M^{2}-2(|a_{1,1}|^{2}+|b_{1,1}|^{2})}=
\sqrt{2M^{2}-2(1+2|a_{1,1}|^{2})}\Big(|a_{1,1}|+\sqrt{1+|a_{1,1}|^{2}}\Big)\cdot \lambda_{F}(0).$$
As in the previous case, we get that
$$\sqrt{2M^{2}-2(1+2|a_{1,1}|^{2})}\Big(|a_{1,1}|+\sqrt{1+|a_{1,1}|^{2}}\Big)\leq\sqrt{M^{4}-1},$$
which proves that the inequality \eqref{eq1.4} holds,
$|a_{n,k}|+|b_{n.k}|\leq \sqrt{M^{4}-1}\cdot \lambda_{F}(0)$ for $(n,k)\not=(1,1)$ and $$\lambda_{F}(0)=\big||a_{1,1}|-|b_{1,1}|\big|=\frac{1}{\sqrt{1+|a_{1,1}|^{2}}+|a_{1,1}|}
\geq\frac{\sqrt{2}}{\sqrt{M^{2}-1}+\sqrt{M^{2}+1}}.$$

Next, we prove the inequality \eqref{eq1.6}. We have
from $\lambda_{F}(0)=\big||a_{1,1}|-|b_{1,1}|\big|=1$
that $|a_{1,1}|^{2}+|b_{1,1}|^{2}\geq1$.
Hence the inequality \eqref{eq1.6} follows from \eqref{eq1.5}.
\epf

From \eqref{eq1.5}, we obtain the following two corollaries.

\begin{cor}\label{cor1}
Suppose that $F$ is a polyharmonic mapping of the form
\eqref{eq1.1} and all its non-zero coefficients $a_{n,k_{1}}$, $a_{n,k_{2}}$ and
$b_{n,k_{3}}$, $b_{n,k_{4}}$ satisfy \eqref{eq1.2}.

If $F$ satisfies
$|F(z)|\leq M$ for all $z\in \mathbb{D}$, then for all $n\in N^{\ast}$,
$$\sum_{k=1}^{p}|a_{n,k}|,\;
\sum_{k=1}^{p}|b_{n,k}|\leq \sqrt{p}M\;\text{and}\;\sum_{k=1}^{p}(|a_{n,k}|+|b_{n,k}|)\leq \sqrt{2p}M.$$

In addition, if $F(0)=0$, $\lambda_{F}(0)=1$ or $|J_{F}(0)|=1$, then $|a_{n,k}|$,
$|b_{n,k}|\leq \sqrt{M^{2}-1}$ and for $n\in N^{\ast} \setminus \{1\}$
 $$\sum_{k=1}^{p}|a_{n,k}|,\;
 \sum_{k=1}^{p}|b_{n,k}|\leq\sqrt{p(M^{2}-1)}\;\text{and}\;\sum_{k=1}^{p}(|a_{n,k}|+|b_{n,k}|)\leq \sqrt{2p(M^{2}-1)}.$$
\end{cor}

\begin{cor}\label{cor2}
Suppose that $F$ is a polyharmonic mapping of the form
\eqref{eq1.1}, all its non-zero coefficients $a_{n,k_{1}}$, $a_{n,k_{2}}$ and
$b_{n,k_{3}}$, $b_{n,k_{4}}$ satisfy \eqref{eq1.2}, $F(0)=0$, $|F(z)|\leq 1$ on $\mathbb{D}$.

 $(I)$ If $J_{F}(0)=1$, then $F(z)=\alpha z$,  where $|\alpha|=1$.
 And if $J_{F}(0)=-1$, then $F(z)=\beta \overline{z}$, where $|\beta|=1$.

 $(II)$ If $\lambda_{F}(0)=1$, then $F(z)=\gamma z$ or $F(z)=\gamma \overline{z}$, where $|\gamma|=1$.

\end{cor}

%

\begin{figure}
\begin{center}

\includegraphics[width=6cm]{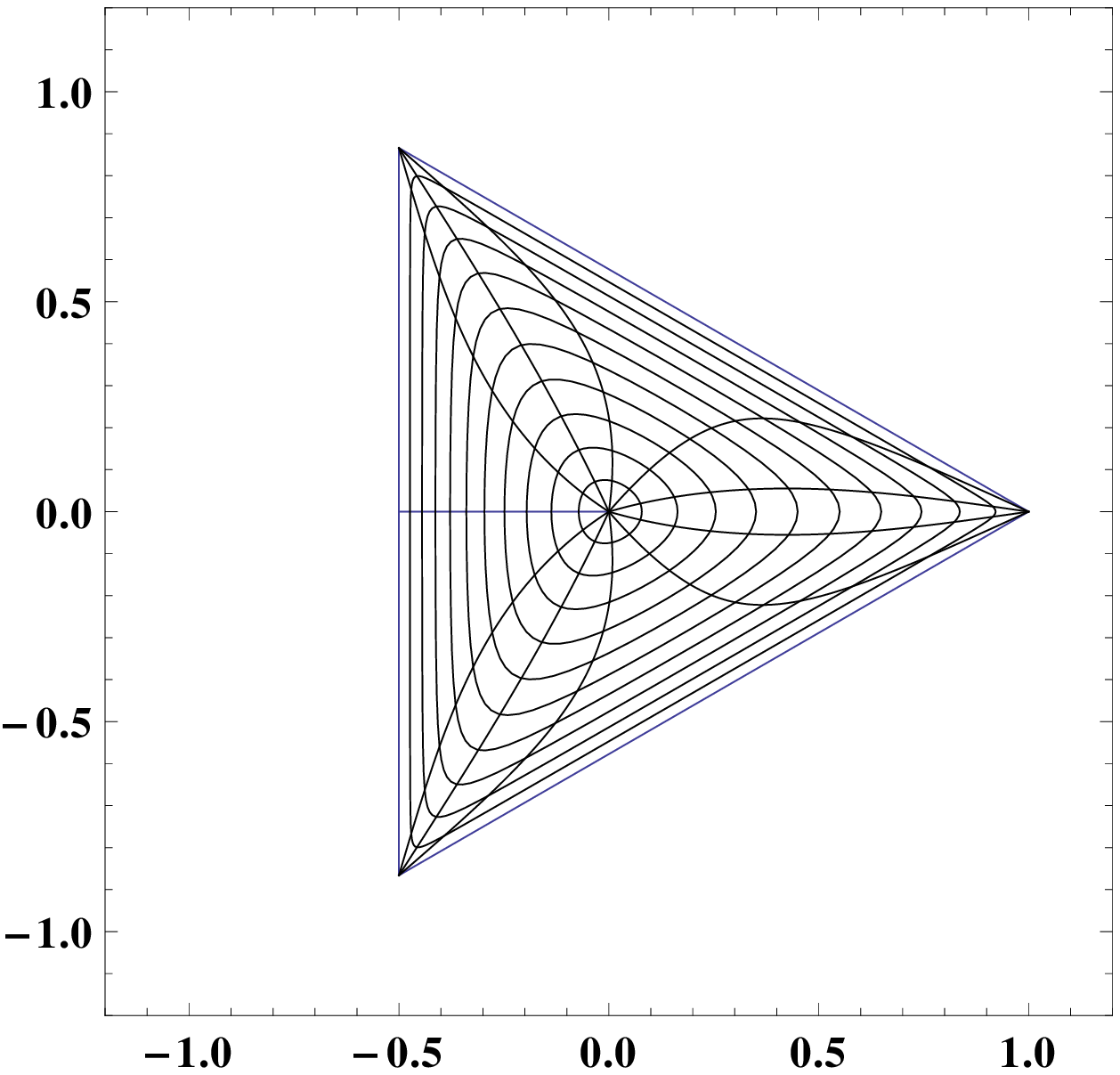}
\quad
\includegraphics[width=6cm]{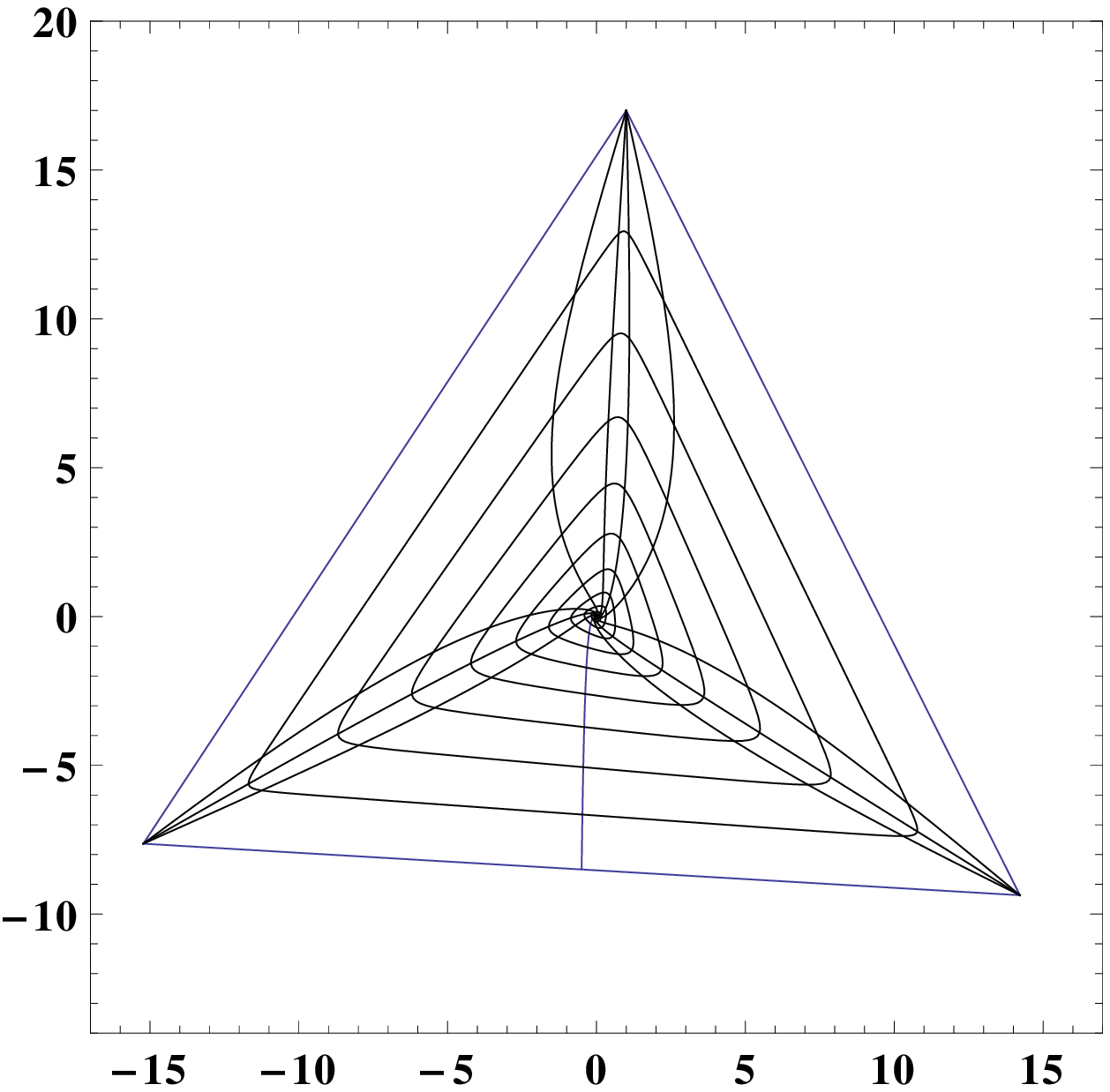}
\end{center}

\caption{The images of $\mathbb{D}$ under $f_{3}$ (left) and $F_{0}$ (right). See Example \ref{example1}.}\label{f3F0}
\end{figure}

\begin{examp}\label{example1} Let  $\beta=\sqrt{\alpha}=e^{i \pi/n}$ and
\begin{eqnarray*}
f_{n}(z)&=&\frac{1}{2\pi i}\sum_{k=0}^{n-1}\alpha^{k}
\left\{\log \frac{z-\beta^{2k+1}}{z-\beta^{2k-1}}-
\overline{\log \frac{z-\beta^{2k+1}}{z-\beta^{2k-1}}}\right\}\\
&=&\sum_{m=1}^{\infty}\big(a_{m}z^{m}+\overline{b_{m}}\overline{z^{m}}\big),\\
\end{eqnarray*}
where $$a_{m}=\begin{cases}
\displaystyle\frac{n}{\pi m}\sin \frac{\pi m}{n},\;\;\;m=1,\;n+1,\;2n+1,\cdots,\\
\displaystyle  0,\;\;\;\;\;\;\;\;\;\;\;\;\;\;\;\;\;\; \text{otherwise},
\end{cases}
$$
$$b_{m}=\begin{cases}
\displaystyle\frac{n}{\pi m}\sin \frac{\pi m}{n},\;\;\;m=n-1,\;2n-1, \cdots,\\
\displaystyle  0,\;\;\;\;\;\;\;\;\;\;\;\;\;\;\;\; \;\;\;\text{otherwise},
\end{cases}
$$ for $n=3,4,\cdots$. Then $f_{n}$ is univalent in $\mathbb{D}$, and maps $\mathbb{D}$ onto
the region inside the regular $n$-gon with vertices 1, $\alpha,\cdots,\alpha^{n-1}$ (cf. \cite[p. 59]{du}).
Let $F_{0}(z)=f_{3}(z)+17i|z|^{2}f_{3}(z)$.
Then $|F_{0}|<18$ in $\mathbb{D}$, and $F_{0}$ satisfies the condition \eqref{eq1.2}. Therefore, by \eqref{eq1.3} we obtain
$$ \sum_{k=1}^{p}\sum_{n=1}^{\infty}(|a_{n,k}|^{2}+|b_{n.k}|^{2})\leq 324,$$
and $|a_{n,k}|+|b_{n,k}|\leq18\sqrt{2} $ for all coefficients.
See Figure \ref{f3F0} for the images of $\mathbb{D}$ under $f_{3}$ and $F_{0}$.
\end{examp}

Now, we establish a version of Landau's theorem for polyharmonic mappings.

\begin{thm}\label{thm2.1}
Suppose that $F$ is a polyharmonic mapping of the form
\eqref{eq1.1} and all its non-zero coefficients $a_{n,k_{1}}$, $a_{n,k_{2}}$ and
$b_{n,k_{3}}$, $b_{n,k_{4}}$ satisfy \eqref{eq1.2}.

If $|F(z)|\leq M$ in $\mathbb{D}$ for some $M>1$
and $F(0)=|J_{F}(0)|-1=0$, then
$F$ is univalent in the disk $\mathbb{D}_{r_{1}}$ and $F(\mathbb{D}_{r_{1}})$
contains a univalent disk $\mathbb{D}_{\rho_{1}}$, where $r_{1}(M,p)$ is
the least positive root of the following equation:
$$1-\sqrt{M^{4}-1}\left(\frac{2r-r^{2}}{(1-r)^{2}}+
\sum_{k=1}^{p-1}\frac{r^{2k}}{(1-r)^{2}} +2 \sum_{k=1}^{p-1}\frac{kr^{2k}}{1-r} \right)=0,$$
and
$$\rho_{1}=\lambda_{0}(M)r_{1}\left(1-\sqrt{M^{4}-1}\frac{r_{1}}{1-r_{1}}-
\sqrt{M^{4}-1}\sum_{k=1}^{p-1}\frac{2r_{1}^{2k}}{1-r_{1}}\right).$$

In particular, if $|F(z)|\leq 1$ for all $z\in \mathbb{D}$ and $F(0)=|J_{F}(0)|-1=0$, then
$F$ is univalent and maps $\mathbb{D}$ onto $\mathbb{D}$.
\end{thm}

\bpf Since $F_{z}(z)$ and $F_{\overline{z}}(z)$
can be written as $$F_{z}(z)=\sum_{k=1}^{p}|z|^{2(k-1)}h_{k}'(z)
+\sum_{k=2}^{p}(k-1)\overline{z}|z|^{2(k-2)}\big(h_{k}(z)+\overline{g_{k}(z)}\big),$$
$$F_{\overline{z}}(z)=\sum_{k=1}^{p}|z|^{2(k-1)}\overline{g_{k}'(z)}
+\sum_{k=2}^{p}(k-1)z|z|^{2(k-2)}\big(h_{k}(z)+\overline{g_{k}(z)}\big),$$
 then for any $z_{1}\not=z_{2}$, where $z_{1}$, $z_{2}\in \mathbb{D}_{r}$ and $r\in (0,1)$ is a constant, we have
\begin{eqnarray*}
\big |F(z_{1})-F(z_{2})\big |&=&\left |\int_{[z_{1},z_{2}]}F_{z}(z)dz+F_{\overline{z}}(z)d\overline{z}
\right |\\
&\geq &\left |\int_{[z_{1},z_{2}]}F_{z}(0)dz+F_{\overline{z}}(0)d\overline{z}\right |\\
&&- \left |\int_{[z_{1},z_{2}]}\big(F_{z}(z)-F_{z}(0)\big)dz
+\big(F_{\overline{z}}(z)-F_{\overline{z}}(0)\big)d\overline{z}
\right |\\
&\geq& J_1-J_2-J_3-J_4,
\end{eqnarray*}
where
\begin{eqnarray*}
J_1& :=& \left |\int_{[z_{1},z_{2}]}h_{1}'(0)\,dz+
\overline{g_{1}'(0)}d\overline{z}\right |,\\
J_2&  :=& \left |\int_{[z_{1},z_{2}]} \big(h_{1}'(z)-h_{1}'(0)\big )dz
+\big (\overline{g_{1}'(z)}- \overline{g_{1}'(0)}\big )d\overline{z}\right|,\\
J_3&  :=& \left |\int_{[z_{1},z_{2}]}\sum_{k=1}^{p-1}|z|^{2k}
 h'_{k+1} (z)dz+\sum_{k=1}^{p-1}|z|^{2k}\overline{ g'_{k+1} (z)}d\overline{z}\right|,\\
J_4& :=& \left |\int_{[z_{1},z_{2}]}\sum_{k=1}^{p-1}k|z|^{2(k-1)}
\big(h_{k+1}(z)+\overline{g_{k+1}(z)}\big)(\overline{z}dz+zd\overline{z})\right|.\\
\end{eqnarray*}
First, we have the estimate for $J_1$,
$$J_1\geq \int_{[z_{1},z_{2}]}\lambda_{F}(0)|dz|
=\lambda_{F}(0)|z_{1}-z_{2}|.
$$
Next, by Theorem \ref{thm1.1} $(II)$, we obtain $|a_{n,k}|+|b_{n,k}|
\leq \sqrt{M^{4}-1}\cdot\lambda_{F}(0)$ for $(n,k)\not=(1,1)$. It follows that
\begin{eqnarray*}
J_2 &\leq&\int_{[z_{1},z_{2}]}\big|h_{1}'(z)-h_{1}'(0)\big | |dz|
+\big |\overline{g_{1}'(z)}- \overline{g_{1}'(0)}\big| |d\overline{z}|\\
&\leq&|z_{1}-z_{2}|\sum_{n=2}^{\infty}n(|a_{n,1}|+|b_{n,1}|)r^{n-1}\\
&\leq&|z_{1}-z_{2}|\sqrt{M^{4}-1}\cdot \lambda_{F}(0)\frac{2r-r^{2}}{(1-r)^{2}},
\end{eqnarray*}
and
\begin{eqnarray*}
J_3 &\leq&\int_{[z_{1},z_{2}]}\sum_{k=1}^{p-1}|z|^{2k}
\big|h_{k+1}'(z)\big| |dz|+\sum_{k=1}^{p-1}|z|^{2k} \big|\overline{g_{k+1}'(z)}\big| |d\overline{z}|\\
&\leq&|z_{1}-z_{2}|\sum_{k=1}^{p-1}r^{2k}\sum_{n=1}^{\infty}n(|a_{n,k+1}|+|b_{n,k+1}|)r^{n-1}\\
&\leq&|z_{1}-z_{2}|\sqrt{M^{4}-1}\cdot \lambda_{F}(0)\sum_{k=1}^{p-1}\frac{r^{2k}}{(1-r)^{2}}.
\end{eqnarray*}
Finally,
\begin{eqnarray*}
J_4&=& \left |\int_{[z_{1},z_{2}]}\sum_{k=1}^{p-1}k|z|^{2(k-1)}
\big(h_{k+1}(z)+\overline{g_{k+1}(z)}\big)(\overline{z}dz+zd\overline{z})\right|\\
&\leq&\left |\int_{[z_{1},z_{2}]}2\sum_{k=1}^{p-1}kr^{2(k-1)}
\sum_{n=1}^{\infty}\left(|a_{n,k+1}|+|b_{n,k+1}|\right)r^{n+1}|dz|\right|\\
&\leq&2|z_{1}-z_{2}|\sum_{k=1}^{p-1}kr^{2(k-1)}
\sum_{n=1}^{\infty}\left(|a_{n,k+1}|+|b_{n,k+1}|\right)r^{n+1}\\
&\leq&2|z_{1}-z_{2}|\sqrt{M^{4}-1}\cdot \lambda_{F}(0)\sum_{k=1}^{p-1}\frac{kr^{2k}}{1-r}.
\end{eqnarray*}
Using these estimates, we obtain
$$  \big |F(z_{1})-F(z_{2})\big | \geq J_1-J_2-J_3-J_4\geq \lambda_{F}(0)|z_{1}-z_{2}|\psi(r),$$
where
\begin{eqnarray*}
\psi(r)&=&1-\sqrt{M^{4}-1}\left(\frac{2r-r^{2}}{(1-r)^{2}}+
\sum_{k=1}^{p-1}\frac{r^{2k}}{(1-r)^{2}} +2 \sum_{k=1}^{p-1}\frac{kr^{2k}}{1-r} \right).
\end{eqnarray*}
It is easy to see that the function $\psi(r)$ is strictly decreasing
for $r\in(0,1)$ and $M>1$,
$$\lim_{r\rightarrow 0+}\psi(r)=1\ \mbox{and}\
\lim_{r\rightarrow 1^{-}}\psi(r)=-\infty.
$$
Hence there exists a unique $r_{1}\in (0,1)$ satisfying
$\psi(r_{1})=0.$ This implies that $F $ is univalent in
$\mathbb{D}_{r_{1}}$.

By Theorem \ref{thm1.1} $(II)$, we see $|a_{n,k}|+|b_{n,k}|\leq
\sqrt{M^{4}-1}\cdot \lambda_{F}(0)$ for all $(n,k)\not=(1,1)$.
Therefore, by \eqref{eq1.7}, for any $w$ in $\{w:\; |w|=r_{1}\}$, we obtain
\begin{eqnarray*}
\big |F(w)-F(0)\big|
&=&\left |\int_{[0,w]}F_{z}(z)dz+F_{\overline{z}}(z)d\overline{z}
\right |\\
&\geq&\lambda_{F}(0)r_{1}-\sum_{n=2}^{\infty}(|a_{n,1}|+|b_{n,1}|)r_{1}^{n}-
2\sum_{k=1}^{p-1}r_{1}^{2k}\sum_{n=1}^{\infty}(|a_{n,k+1}|+|b_{n,k+1}|)r_{1}^{n}\\
&\geq&\lambda_{F}(0)r_{1}\left(1-\sqrt{M^{4}-1}\frac{r_{1}}{1-r_{1}}-
\sqrt{M^{4}-1}\sum_{k=1}^{p-1}\frac{2r_{1}^{2k}}{1-r_{1}}\right)\\
&\geq&\lambda_{0}(M)r_{1}\left(1-\sqrt{M^{4}-1}\frac{r_{1}}{1-r_{1}}-
\sqrt{M^{4}-1}\sum_{k=1}^{p-1}\frac{2r_{1}^{2k}}{1-r_{1}}\right):=\rho_{1}.\\
\end{eqnarray*}
We also have $$\rho_{1}>\lambda_{0}(M)r_{1}\left[1-\sqrt{M^{4}-1}\left(\frac{2r_{1}-r_{1}^{2}}{(1-r_{1})^{2}}+
\sum_{k=1}^{p-1}\frac{r_{1}^{2k}}{(1-r_{1})^{2}} +2 \sum_{k=1}^{p-1}\frac{kr_{1}^{2k}}{1-r_{1}} \right)\right]=0.$$

If $|F(z)|\leq 1$ for all $z\in \mathbb{D}$ and $F(0)=|J_{F}(0)|-1=0$,
then by Corollary \ref{cor2} ($I$), we have $F$ is univalent and maps $\mathbb{D}$ onto $\mathbb{D}$.
The proof of the theorem is complete.
\epf

\begin{cor}\label{cor2.2}
Suppose that $F$ is a polyharmonic mapping of the form
\eqref{eq1.1} and all its non-zero coefficients $a_{n,k_{1}}$, $a_{n,k_{2}}$ and
$b_{n,k_{3}}$, $b_{n,k_{4}}$ satisfy \eqref{eq1.2}.

If $|F(z)|\leq M$ in $\mathbb{D}$ for some $M>1$
and $F(0)=\lambda_{F}(0)-1=0$, then
$F$ is univalent in the disk $\mathbb{D}_{r_{2}}$ and $F(\mathbb{D}_{r_{2}})$
contains a univalent disk $\mathbb{D}_{\rho_{2}}$, where $r_{2}(M,p)$ is
the least positive root of the following equation:
$$1-\sqrt{2M^{2}-2}\left(\frac{2r-r^{2}}{(1-r)^{2}}+
\sum_{k=1}^{p-1}\frac{r^{2k}}{(1-r)^{2}} +2 \sum_{k=1}^{p-1}\frac{kr^{2k}}{1-r} \right)=0,$$
and
$$\rho_{2}=r_{2}\left(1-\sqrt{2M^{2}-2}\frac{r_{2}}{1-r_{2}}-
\sqrt{2M^{2}-2}\sum_{k=1}^{p-1}\frac{2r_{2}^{2k}}{1-r_{2}}\right).$$

In particular, if $|F(z)|\leq 1$ for all $z\in \mathbb{D}$ and $F(0)=\lambda_{F}(0)-1=0$, then
$F$ is univalent and maps $\mathbb{D}$ onto $\mathbb{D}$.
\end{cor}
\bpf The proof of this result is similar to Theorem \ref{thm1.1}, where
$|a_{n,k}|+|b_{n,k}|\leq \sqrt{2M^{2}-2}$ and $\lambda_{F}(0)=1$ is used
instead of $|a_{n,k}|+|b_{n,k}|\leq \sqrt{M^{4}-1}\cdot \lambda_{F}(0)$
for all $(n,k)\not=(1,1)$, and we omit it.
\epf

\begin{examp} Let $F_{1}(z)=K(z)+|z^{2}|G(z)=\frac{2\pi}{3\sqrt{3}}f_{3}(z)+\frac{34\pi i}{3\sqrt{3}}|z|^{2}f_{3}(z)$ 
for $z\in \mathbb{D} $.
Then $F_{1}(0)=J_{F_{1}}(0)-1=\lambda_{F_{1}}(0)-1=0$, $|F_{1}(z)|<M_{1}:=4\sqrt{3}\pi$ and
both $|K|$ and $|G|$ are bounded by $M_{2}:=\frac{34\pi }{3\sqrt{3}}$
. Then, by Corollary \ref{cor2.2},
$F_{1}$ is univalent in the disk $\mathbb{D}_{r_{3}}$ and $F(\mathbb{D}_{r_{3}})$
contains a univalent disk $\mathbb{D}_{\rho_{3}}$, where $r_{3}\approx0.01552$ is
the least positive root of the following equation:
$$1-\sqrt{2M_{1}^{2}-2}\left(\frac{2r }{(1-r)^{2}}+2 \frac{ r^{2 }}{1-r} \right)=0,$$
and
$$\rho_{3}= r_{3} \left(1-\sqrt{2M_{1}^{2}-2}\cdot\frac{r_{3}+2r_{3} ^{2 }}{1-r_{3} }\right)\approx0.00776,$$
while by \cite[Theorem 2]{sh2011}, we see that
$F_{1}$ is univalent in the disk $\mathbb{D}_{r_{4}}$ and $F_{1}(\mathbb{D}_{r_{4}})$
contains a univalent disk $\mathbb{D}_{\rho_{4}}$, where $r_{4}\approx0.00041$ is
the least positive root of the following equation:
$$\frac{\pi}{4M_{2}}-\frac{4M_{2}r(2-r)}{\pi(1-r)^{2}}-\frac{4M_{2}r^{2}}{\pi(1-r)^{2}}-2M_{2}r=0,$$
and
$$\rho_{4}= r_{4} \left(
\frac{\pi}{4M_{2}}-\frac{4M_{2}r_{4}}{\pi (1-r_{4})}-\frac{4M_{2}r_{4}^{2}}{\pi(1-r_{4})}
\right)\approx 1.12385\times 10^{-5}.$$
\end{examp}

\begin{cor}\label{cor2.1}
Suppose that $F$ is a polyharmonic mapping of the form
$$F(z)=\sum_{k=1}^{p}\lambda_{k}|z|^{2(k-1)}G(z)=\sum_{k=1}^{p}|z|^{2(k-1)}
\sum_{n=1}^{\infty}(a_{n}z^{n}+\overline{b_{n}}\overline{z^{n}}),$$
where $G(z)$ is nonconstant, $\sum_{k=1}^{p}\lambda_{k}=1$ and  $\lambda_{k}\in [0,1]$ for all $k\in\{1,\cdots,p\}$.
If $|F(z)|\leq M$ in $\mathbb{D}$ for some $M>1$ and $\lambda_{F}(0)=1$, then
$F$ is univalent in the disk $\mathbb{D}_{r_{5}}$, and $F(\mathbb{D}_{r_{5}})$
contains a univalent disk $\mathbb{D}_{\rho_{5}}$, where $r_{5}(M,p)$ is
the least positive root of the following equation:
$$1-T(M)\left(\frac{2r-r^{2}}{(1-r)^{2}}+
\sum_{k=1}^{p-1}\frac{r^{2k}}{(1-r)^{2}} +2 \sum_{k=1}^{p-1}\frac{kr^{2k}}{1-r} \right)=0,$$
and
$$\rho_{5}=r_{5}\left(1-T(M)\frac{r_{5}}{1-r_{5}}-
T(M)\sum_{k=1}^{p-1}\frac{2r_{5}^{2k}}{1-r_{5}}\right),$$
where $T(M):=\min \left\{\sqrt{2M^{2}-2}, \;\frac{4M}{\pi} \right\}$.
 
\end{cor}
\bpf Suppose that $|F(z)|\leq M$ for all $z\in \mathbb{D}$ and $\sum_{k=1}^{p}\lambda_{k}=1$.
Then $$\lim_{r\rightarrow 1^{-}}|F(re^{i\theta})|=\lim_{r\rightarrow 1^{-}}|G(re^{i\theta})|\leq M.$$
Since $G(z)$ is harmonic and nonconstant, by the maximum principle,
we have $|G(z)|<M$ in $\mathbb{D}$.
Lemma B and Theorem \ref{thm1.1} ($III$) implies that
$$|a_{n}|+|b_{n}|\leq T(M)= \min\left\{\sqrt{2M^{2}-2},\frac{4}{\pi}M\right\}$$
for all $n\geq1$. By Corollary \ref{cor2.2}, we easily obtain the results.
\epf

Now, we consider Landau's theorem for the mapping $L(F)$, where $F$ is polyharmonic.
\begin{thm}\label{thm3.1}
Suppose that $F$ is a polyharmonic mapping of the form
\eqref{eq1.1} and all its non-zero coefficients $a_{n,k_{1}}$, $a_{n,k_{2}}$ and
$b_{n,k_{3}}$, $b_{n,k_{4}}$ satisfy \eqref{eq1.2}.

If $|F(z)|\leq M$ in $\mathbb{D}$ for some $M>1$
and $F(0)=|J_{F}(0)|-1=0$, then
$L(F)$ is univalent in the disk $\mathbb{D}_{r_{6}}$ and $F(\mathbb{D}_{r_{6}})$
contains a univalent disk $\mathbb{D}_{\rho_{6}}$, where $r_{6}(M,p)$ is
the least positive root of the following equation:
$$1-\sqrt{M^{4}-1}\left(\frac{2r-r^{2}}{(1-r)^{2}}+
\sum_{k=1}^{p}\frac{2r^{2k-1}}{(1-r)^{3}} +\sum_{k=2}^{p}\frac{(2k-1)r^{2(k-1)}}{(1-r)^{2}} \right)=0,$$
and
$$\rho_{6}=\lambda_{0}(M)r_{6}\left(1-\sqrt{M^{4}-1}\frac{2r_{6}-r_{6}^{2}}{(1-r_{6})^{2}}-
\sqrt{M^{4}-1}\sum_{k=2}^{p}\frac{r_{6}^{2(k-1)}}{(1-r_{6})^{2}}\right).$$

In particular, if $|F(z)|\leq 1$ for all $z\in \mathbb{D}$ and $F(0)=|J_{F}(0)|-1=0$, then
$L(F)$ is univalent and maps $\mathbb{D}$ onto $\mathbb{D}$.
\end{thm}

\bpf Let $$H(z):=L(F)=zF_{z}(z)-\overline{z}F_{\overline{z}}(z)
=\sum_{k=1}^{p}|z|^{2(k-1)}\Big(zh'_{k}(z)-\overline{z}\overline{g'_{k} (z)}\Big).$$
Then $H_{z}(z)$ and $H_{\overline{z}}(z)$ can be written as
\begin{eqnarray*}
H_{z}(z)&=&h'_{1} (z)+z h''_{1} (z)
+\sum_{k=2}^{p}\Big(k|z|^{2(k-1)}h'_{k}(z)+z|z|^{2(k-1)}h''_{k}(z)\\
&&-(k-1)\overline{z^{2}}|z|^{2(k-2)}\overline{g'_{k}(z)}\Big),\\
\end{eqnarray*}
and\begin{eqnarray*}
-H_{\overline{z}}(z)&=&\overline{g'_{1}(z)}+\overline{z}\overline{g''_{1}(z)}
+\sum_{k=2}^{p}\Big(k|z|^{2(k-1)}\overline{g'_{k}(z)}+\overline{z}|z|^{2(k-1)}
\overline{g''_{k}(z)}\\
&&-(k-1)z^{2}|z|^{2(k-2)}h'_{k}(z)\Big).\\
\end{eqnarray*}
Let $r$ be a constant in $(0,1)$. For any distinct $z_{1}$, $z_{2}\in \mathbb{D}_{r}$, we have
\begin{eqnarray*}
\big |H(z_{1})-H(z_{2})\big |&=&\left |\int_{[z_{1},z_{2}]}H_{z}(z)dz+H_{\overline{z}}(z)d\overline{z}
\right |\\
&\geq &\left |\int_{[z_{1},z_{2}]}H_{z}(0)dz+H_{\overline{z}}(0)d\overline{z}\right |\\
&&- \left |\int_{[z_{1},z_{2}]}\big(H_{z}(z)-H_{z}(0)\big)dz+\big(H_{\overline{z}}(z)-H_{\overline{z}}(0)\big)d\overline{z}
\right |\\
&\geq& K_1-K_2-K_3-K_4,
\end{eqnarray*}
where
\begin{eqnarray*}
K_1& :=& \left |\int_{[z_{1},z_{2}]}h'_{1}(0)\,dz-
\overline{g'_{1}(0)}d\overline{z}\right |,\\
K_2& :=& \left |\int_{[z_{1},z_{2}]}
\big(h'_{1}(z)-h'_{1}(0)\big)dz-
\big(\overline{g'_{1}(z)}-\overline{g'_{1}(0)}\big)d\overline{z}\right|,\\
K_3& :=& \left |\int_{[z_{1},z_{2}]} \sum_{k=1}^{p}|z|^{2(k-1)}
\Big(z h''_{k}(z) dz
-\overline{z}\overline{g''_{k}(z)}d\overline{z}
\Big )\right|,\\
K_4& :=& \left |\int_{[z_{1},z_{2}]}\sum_{k=2}^{p}(k-1)|z|^{2(k-2)}
\Big(z^{2}h'_{k}(z)d\overline{z}-\overline{z^{2}}\overline{g'_{k}(z)}dz\Big )\right|\\
&&+\left |\int_{[z_{1},z_{2}]}\sum_{k=2}^{p}k|z|^{2(k-1)}
\Big(h'_{k}(z)dz-\overline{g'_{k}(z)}d\overline{z}\Big )\right|.\\
\end{eqnarray*}
First, observe the following estimate for $K_1$:
$$K_1\geq \int_{[z_{1},z_{2}]}\lambda_{F}(0)|dz|
=\lambda_{F}(0)|z_{1}-z_{2}|.
$$
Next, by Theorem \ref{thm1.1} $(II)$, we have
$|a_{n,k}|+|b_{n,k}|\leq \sqrt{M^{4}-1}\cdot \lambda_{F}(0)$ for $(n,k)\not=(1,1)$. Hence
\begin{eqnarray*}
 K_2&=&\left |\int_{[z_{1},z_{2}]} \sum_{n=2}^{\infty}na_{n,1}z^{n-1} dz-
 \sum_{n=2}^{\infty}n\overline{b_{n,1}} \overline{z^{n-1}} d\overline{z}  \right|\\
&\leq&|z_{1}-z_{2}|\sum_{n=2}^{\infty}n(|a_{n,1}|+|b_{n,1}|)r^{n-1}\\
&\leq&|z_{1}-z_{2}|\sqrt{M^{4}-1}\cdot \lambda_{F}(0)\frac{2r-r^{2}}{(1-r)^{2}},\\
\end{eqnarray*}and
\begin{eqnarray*}
K_3 &\leq&|z_{1}-z_{2}|\sum_{k=1}^{p}r^{2(k-1)}\sum_{n=2}^{\infty}n(n-1)(|a_{n,k}|+|b_{n,k}|)r^{n-1}\\
&\leq&|z_{1}-z_{2}|\sqrt{M^{4}-1}\cdot \lambda_{F}(0)\sum_{k=1}^{p}\frac{2r^{2k-1}}{(1-r)^{3}}.
\end{eqnarray*}
Finally,
\begin{eqnarray*}
K_4 &\leq&\int_{[z_{1},z_{2}]}\sum_{k=2}^{p}(k-1)r^{2(k-2)}\sum_{n=1}^{\infty}n(|a_{n,k}|+|b_{n,k}|)r^{n+1}|dz|\\
&&+\int_{[z_{1},z_{2}]}\sum_{k=2}^{p}kr^{2(k-1)}\sum_{n=1}^{\infty}n(|a_{n,k}|+|b_{n,k}|)r^{n-1}|dz|\\
&\leq&|z_{1}-z_{2}|\sum_{k=2}^{p}(2k-1)r^{2(k-1)}\sum_{n=1}^{\infty}n(|a_{n,k}|+|b_{n,k}|)r^{n-1}\\
&\leq&|z_{1}-z_{2}|\sqrt{M^{4}-1}\cdot \lambda_{F}(0)\sum_{k=2}^{p}(2k-1)\frac{r^{2(k-1)}}{(1-r)^{2}}.
\end{eqnarray*}
Using these estimates, we obtain
$$  \big |H(z_{1})-H(z_{2})\big | \geq K_1-K_2-K_3-K_4\geq \lambda_{F}(0)|z_{1}-z_{2}|\varphi(r),$$
where
\begin{eqnarray*}
\varphi(r)&=&1-\sqrt{M^{4}-1}\left(\frac{2r-r^{2}}{(1-r)^{2}}+
\sum_{k=1}^{p}\frac{2r^{2k-1}}{(1-r)^{3}} +\sum_{k=2}^{p}\frac{(2k-1)r^{2(k-1)}}{(1-r)^{2}} \right).
\end{eqnarray*}
It is easy to verify that $\varphi(r)$ is a strictly decreasing function for $r\in(0,1)$ and $M>1$,
$$\lim_{r\rightarrow 0+}\varphi(r)=1\ \mbox{and}\
\lim_{r\rightarrow 1^{-}}\varphi(r)=-\infty.
$$
Hence there exists a unique $r_{6}\in (0,1)$ satisfying
$\varphi(r_{6})=0.$ This shows that $H(z)$ is univalent in
$\mathbb{D}_{r_{6}}$.

Furthermore, for any $w$ in $\{w:\; |w|=r_{6}\}$,
\begin{eqnarray*}
\big |H(w)\big|
&=&\left|\sum_{k=1}^{p}|w|^{2(k-1)}\big(w h'_{k}(w)-\overline{w}\overline{g'_{k}(z)}\big) \right|\\
&\geq&\big|wh'_{1}(0)-\overline{w}\overline{g'_{1}(0)}\big|-
\Big|w\big( h'_{1}(w)-  h'_{1}(0)\big)-\overline{w}\big( \overline{g'_{1}(w)}-\overline{g'_{1}(0)}\big)\Big|\\
&&-\left|\sum_{k=2}^{p}|w|^{2(k-1)}\big( wh'_{k}(w)-\overline{w}\overline{ g'_{k} (w)} \big)  \right|\\
&\geq&\lambda_{F}(0)r_{6}-\sum_{n=2}^{\infty}n(|a_{n,1}|+|b_{n,1}|)r_{6}^{n}-\sum_{k=2}^{p}r_{6}^{2(k-1)}
\sum_{n=1}^{\infty}n(|a_{n,k}|+|b_{n,k}|)r_{6}^{n}.\\
\end{eqnarray*}
By Theorem \ref{thm1.1} $(II)$, we see $|a_{n,k}|+|b_{n,k}|\leq
\sqrt{M^{4}-1}\cdot \lambda_{F}(0)$ for all $(n,k)\not=(1,1)$.
Therefore, we obtain
\begin{eqnarray*}
\big |H(w)\big|&\geq&\lambda_{F}(0)r_{6}\left(1-\sqrt{M^{4}-1}\frac{2r_{6}-r_{6}^{2}}{(1-r_{6})^{2}}-
\sqrt{M^{4}-1}\sum_{k=2}^{p}\frac{r_{6}^{2(k-1)}}{(1-r_{6})^{2}}\right)\\
&\geq&\lambda_{0}(M)r_{6}\left(1-\sqrt{M^{4}-1}\frac{2r_{6}-r_{6}^{2}}{(1-r_{6})^{2}}-
\sqrt{M^{4}-1}\sum_{k=2}^{p}\frac{r_{6}^{2(k-1)}}{(1-r_{6})^{2}}\right):=\rho_{6},\\
\end{eqnarray*}
where $\lambda_{0}(M)$ is the same as in \eqref{eq1.7}, and
\begin{eqnarray*}
\rho_{6}&>&\lambda_{0}(M)r_{6}\left[1-\sqrt{M^{4}-1}\left(\frac{2r_{6}-r_{6}^{2}}{(1-r_{6})^{2}}+\sum_{k=1}^{p}\frac{2r_{6}^{2k-1}}{(1-r_{6})^{3}} \right.\right.\\
&&\left.\left.+\sum_{k=2}^{p}\frac{(2k-1)r_{6}^{2(k-1)}}{(1-r_{6})^{2}} \right)\right]=0.
\end{eqnarray*}

If $|F(z)|\leq 1$ for all $z\in \mathbb{D}$ and $F(0)=|J_{F}(0)|-1=0$,
then by Corollary \ref{cor2} ($I$), we have $L(F)$ is univalent and maps $\mathbb{D}$ onto $\mathbb{D}$.
The proof of the theorem is complete.
\epf

\begin{cor}\label{cor3.2}
Suppose that $F$ is a polyharmonic mapping of the form
\eqref{eq1.1} and all its non-zero coefficients $a_{n,k_{1}}$, $a_{n,k_{2}}$ and
$b_{n,k_{3}}$, $b_{n,k_{4}}$ satisfy \eqref{eq1.2}.

If $|F(z)|\leq M$ in $\mathbb{D}$ for some $M>1$
and $F(0)=\lambda_{F}(0)-1=0$, then
$L(F)$ is univalent in the disk $\mathbb{D}_{r_{7}}$ and $F(\mathbb{D}_{r_{7}})$
contains a univalent disk $\mathbb{D}_{\rho_{7}}$, where $r_{7}(M,p)$ is
the least positive root of the following equation:
$$1-\sqrt{2M^{2}-2}\left(\frac{2r-r^{2}}{(1-r)^{2}}+
\sum_{k=1}^{p}\frac{2r^{2k-1}}{(1-r)^{3}} +\sum_{k=2}^{p}\frac{(2k-1)r^{2(k-1)}}{(1-r)^{2}} \right)=0,$$
and
$$\rho_{7}=r_{7}\left(1-\sqrt{2M^{2}-2}\frac{2r_{7}-r_{7}^{2}}{(1-r_{7})^{2}}-
\sqrt{2M^{2}-2}\sum_{k=2}^{p}\frac{r_{7}^{2(k-1)}}{(1-r_{7})^{2}}\right).$$

In particular, if $|F(z)|\leq 1$ for all $z\in \mathbb{D}$ and $F(0)=\lambda_{F}(0)-1=0$, then
$L(F)$ is univalent and maps $\mathbb{D}$ onto $\mathbb{D}$.
\end{cor}

\begin{examp} Under the hypothesis of Example 2, and by Corollary \ref{cor3.2},
$L(F_{1})$ is univalent in the disk $\mathbb{D}_{r_{8}}$ and $F_{1}(\mathbb{D}_{r_{8}})$
contains a univalent disk $\mathbb{D}_{\rho_{8}}$, where $r_{8}\approx0.00798$ is
the least positive root of the following equation:
$$1-\sqrt{2M_{1}^{2}-2} \cdot\frac{4r-3r^{2}+3r^{3}+3r^{4}-3r^{5} }{(1-r)^{3}}=0,$$
and
$$\rho_{8}= r_{8} \left(1-\sqrt{2M_{1}^{2}-2}\cdot\frac{2r_{8}}{(1-r_{8})^{2} }\right)\approx0.00400,$$
while by \cite[Theorem 1.1]{sh2009}, we see that
$L(F_{1})$ is univalent in the disk $\mathbb{D}_{r_{9}}$ and $F(\mathbb{D}_{r_{9}})$
contains a univalent disk $\mathbb{D}_{\rho_{9}}$, where $r_{9}\approx0.00013$ is
the least positive root of the following equation:
$$\frac{\pi}{4M_{2}}-\frac{6M_{2}r^{2}}{(1-r)^{2}}-\frac{4Mr^{3}}{(1-r)^{3}}
-\frac{16M_{2}}{\pi^{2}}m_{1}\arctan r-\frac{4M_{2}r}{(1-r)^{3}}=0,$$
where $m_{1}\approx 6.05934$ is the minimum value of the function
$$\frac{2-x^{2}+\frac{4}{\pi}\arctan x}{x(1-x^{2})}$$
for $0<x<1$, and
$$\rho_{9}= r_{9} \left(\frac{\pi}{4M_{2}} -\frac{2M_{2}r_{9}^{2}}{(1-r_{9})^{2}}
-\frac{16M_{2}}{\pi^{2}}m_{1}\arctan r_{9}\right)\approx 1.48687\times 10^{-6}
.$$

\end{examp}

\end{document}